\documentclass{amsart}
\usepackage{amssymb}

\newtheorem{theo}{Theorem}
\newtheorem{thm}{Theorem}[section]
\newtheorem{cor}[thm]{Corollary}

\newtheorem{lem}[thm]{Lemma}
\newtheorem{prop}[thm]{Proposition}

{}
\theoremstyle{remark}
\newtheorem{rem}{Remark}

\def \eqskip { \vspace*{3mm} }

\newcommand{\ds}{\displaystyle}
\newcommand{\dsum}{\ds\sum}
\newcommand{\dint}{\ds\int}

\newcommand{\pol}[1]{{\rm Li}_{#1}}

\begin{document}

\title[Integrals of polylogarithmic functions]
{Integrals of polylogarithmic functions, recurrence relations, 
and associated Euler sums}

\author{Pedro Freitas}

\address{Departamento de Matem\'{a}tica, Instituto Superior 
T\'{e}cnico,
Av. Rovisco Pais, 1049-001 Lisboa, Portugal}
\email{pfreitas@math.ist.utl.pt}
\thanks{Partially supported by FCT, Portugal, through program POCTI}



\date{\today}

\dedicatory{}


\begin{abstract}
We show that integrals of the form
\[
\dint_{0}^{1} x^{m}\pol{p}(x)\pol{q}(x)dx, \;\;(m\geq -2, p,q\geq 1)
\]
and
\[
\dint_{0}^{1} \frac{\ds \log^{r}(x)\pol{p}(x)\pol{q}(x)}{\ds x}dx,
\;\;(p,q,r\geq 1)
\]
satisfy certain recurrence relations which allow us to write them in 
terms of Euler sums. From this we prove that, in the first case for all 
$m,p,q$, and in the second when $p+q+r$ is even, these integrals are reducible
to zeta values. In the case of odd $p+q+r$, we combine the known results 
for Euler sums with the information obtained from the problem in this 
form, to give an estimate on the number of {\it new} constants which 
are needed to express the above integrals for a given weight $p+q+r$.

The proofs are constructive, giving a method for the evaluation
of these and other similar integrals, and we present a selection of 
explicit evaluations in the last section.
\end{abstract}

\maketitle

\section{Introduction}

The evaluation of series of the form
\[
S_{r^p,q}:=\dsum_{n=1}^{\infty} \frac{\ds 
\left[H_{n}^{(r)}\right]^{p}}{\ds n^{q}},
\mbox{ where }
H_{n}^{(r)} := \dsum_{k=1}^{n} \frac{\ds 1}{\ds k^{r}},\; \; 
H_{n}\equiv H_{n}^{(1)},
\]
dates back to Euler (see the introductory note in~\cite{bbg}), and 
has recently been the object of active 
research~\cite{baap,bobo,bbg,chu,doel,flsa}. Apart from the actual 
evaluation of the series, one of the main 
questions that one sets out to solve is whether or not a given series 
can be expressed in terms of a linear rational combination of known constants.
When this is the case, we say that the series is reducible to these 
values.

The purpose of this paper is to carry out a similar analysis in 
the case of integrals involving polylogarithmic 
functions and which can be related to Euler sums of the above type.
We then use these relations to establish whether 
or not these integrals are reducible to known constants, and to 
evaluate them.  

Within this scope, a typical result that can be obtained by using
the techniques in this paper is the following
\begin{theo}\label{maint}
    Let $p,q$ and $m$ be integers with $p,q\geq 1$ and 
    $m\geq -2$. Then the integrals
    \[
    J(m,p,q) := \dint_{0}^{1} x^{m}\pol{p}(x)\pol{q}(x)dx
    \]
    are reducible to a rational 
    constant and zeta values at positive integers.
\end{theo}

As usual, we have denoted by $\pol{p}$ the polylogarithmic function 
defined for $|x|<1$ by
\[
\pol{p}(x)=\dsum_{k=1}^{\infty}\frac{\ds x^{k}}{\ds k^{p}}.
\]

A second result is
\begin{theo}\label{maint2}
    Let $p,q$ and $r$ be integers such that $p,q\geq 0$ with at least 
    one of them positive, $r\geq 1$ and $w=p+q+r$ even.
    Then the integrals
    \[
    K(r,p,q) := \dint_{0}^{1} \frac{\ds 
    \log^{r}(x)\pol{p}(x)\pol{q}(x)}{\ds x}dx
    \]
    are reducible to zeta values at positive integers.
    
    In the case of odd weight $w$, these integrals are reducible 
    to zeta values at positive integers and at 
    most $\lfloor (w-1)/6\rfloor$
    other constants which can be taken from the set
    \[
    \mathcal{S}_{w}=\left\{ \kappa_{2j-1,p+q+r-2j+1}, j=1,\ldots,\lfloor 
    (w+1)/4\rfloor\right\},
    \]
    where $\kappa_{rq}:= K(r,0,q)/r!$.
\end{theo}
\begin{rem} {\rm For low odd weights ($w<7$), it is still possible to reduce 
the integrals only to zeta values, since the integrals
mentioned 
above will also be reducible. See the discussion following 
Theorem~\ref{r0q} and also the examples in Section~\ref{examp}, where the 
cases of weight $3$ and $5$ are given. There we also include the case 
of weight $7$ where all integrals can be written in terms of 
$K(1,0,6)$.}
\end{rem}
\begin{rem} {\rm We actually prove that if one considers the set of 
zeta values at positive integers together with the constants in the set
$\mathcal{S}_{w}$
above, then a specific integral $K(r,p,q)$ can be reduced to zeta 
values and at most $\lfloor (r+1)/2\rfloor$ elements of the set 
$\mathcal{S}_{w}$ in the case of $r$ smaller than or equal to $(w-1)/2$, 
and at most $\lfloor (w-r)/2\rfloor$ elements of 
$\mathcal{S}_{w}$ otherwise.
}
\end{rem}

Although in a completely different spirit from the rest of the paper, 
we also give a very simple approximate expression for the elements of the set 
$\mathcal{S}_{w}$. More precisely, we have
\begin{theo}\label{maint3}
    The constants $\kappa_{rq}$ satisfy
    the following estimate
    \[
    \left| \kappa_{rq}-(-1)^{r} \zeta(q)\left[\zeta(r+1)-1\right]\right|\leq
    \frac{\ds \zeta(q-1)-\zeta(q)}{\ds 2^{r+1}}, \;\; r\geq 1\;,q\geq 3.
    \]
\end{theo}

Integrals of the type considered here have been studied in the 
literature either in connection with the general study of
polylogarithms \cite{lewi}, or due to the role they play in calculations
involving higher order Feynman diagrams, 
relativistic phase space integrals for scattering processes under 
some restrictions, etc \cite{dedu,gatr}. In particular, the connection between some of these
integrals and Euler series had already been explored in~\cite{gatr}, although
there the authors restricted their study to integrals containing only
logarithms such as
\[
\int_{0}^{1} \frac{\ds \log^{m}(x)\log^{n}(1-x)\log^{p}(1+x)}{\ds 
a(x)}dx,
\]
where $a(x) = x,1-x,$ or $1+x$, and $m,n$ and $p$ are integers. On 
the other hand, this line of research had also been suggested 
in~\cite{srch}, where the authors point out the relation between 
linear Euler sums and integrals of the type of $K(r,0,q)$ above -- see 
Remark 3 on page 157 of~\cite{srch} and Lemma~\ref{transf2} below.

In the present paper, and although we do give some explicit 
evaluations (see Theorems~\ref{r0q},~\ref{mult} and~\ref{jm1pq},
Corollary~\ref{cor2} and 
Lemma~\ref{symmrec}), the emphasis is on the determination of a set of
recurrence relations that is satisfied by the integrals. We believe that this
point of view has certain advantages, as in many cases of this type the
expressions
for the integrals quickly get quite large, and while getting a 
general closed expression may be quite cumbersome, it is quite simple to 
write down a set of rules that can then be run on a software package such as 
Maple or Mathematica and which will then make it possible to evaluate
any specific integral that is needed. Of course in order to 
proceed in this way, one then needs to carry out the explicit
evaluation of the {\it initial conditions} for these relations. One of the 
ingredients in this last step in some of the cases considered are the 
relations already mentioned
between the integrals and certain Euler sums.

Proceeding in a similar fashion to the evaluation of the integrals $K$ and $J$,
it is possible to evaluate other integrals involving polylogarithms, 
and we shall now mention some further examples. These include, for instance,
certain multiple integrals related to the sums $S_{1^{p},q}$, and integrals of
the form
\[
L(m,r,p) := \dint_{0}^{1} x^{m}\log^{r}(x)\pol{p}(x)dx, \;\; m\geq -1,
\]
for which it is quite straightforward to obtain that
\[
L(-1,r,p) = (-1)^{r}r!\zeta(p+r+1),
\]
and then to use recurrence relations in a similar fashion to 
what was done for $J$ to show that the remaining cases can also be reduced to
zeta values at positive integers plus a rational constant.

The plan of the paper is as follows. In Section~\ref{transfsec} we 
introduce the relations between the integrals and Euler sums. These 
are then used to evaluate some integrals 
directly in Section~\ref{redsec}, where we also establish the recurrence
relations which are needed to prove Theorems~\ref{maint} and~\ref{maint2}.
This is done in Section~\ref{pt12}, and the estimates in 
Theorem~\ref{maint3} are proven in Section~\ref{pt3}.
In the last 
section, 
we present a collection of some explicit evaluations. Finally, in the Appendix
we give an example where the procedure is
reversed, and where we use the corresponding integral to evaluate the Euler sum
$S_{1^{2},2}$.

\section{\label{transfsec}Identities between integrals of polylogarithms and
Euler sums}

We begin by establishing some relations between Euler sums and certain 
integrals. Although it is possible to obtain more general 
identities by proceeding in the same fashion, for simplicity we have chosen
to highlight only some possibilities. The first is a known relation 
between linear Euler sums and $K(r-1,0,q)$ which can be found, for instance, 
in~\cite{srch} (Proposition 3.7). For completeness, and since this 
will play an important role in what follows, we provide a (different)
proof of this result.

\begin{lem}\label{transf2}
    For $q,r\geq 2$ we have that
    \[
    S_{r,q}=\zeta(r)\zeta(q)-\frac{(-1)^{r-1}}{\ds 
    (r-1)!}\int_{0}^{1}\frac{\ds \log^{r-1}(x)\pol{q}(x)}{\ds 1-x}dx.
    \]
\end{lem}
\begin{proof}
    From the fact that
    \begin{equation}\label{xlogint}
    \int_{0}^{1} x^{k-1}\log^{r-1}(x)dx = (-1)^{r-1}\frac{\ds 
    (r-1)!}{k^{r}} \;\;\; (k\geq 1, r\geq 2),
    \end{equation}
    we have that
    \[
    \sum_{k=1}^{n}\frac{\ds 1}{\ds k^{r}} =\frac{\ds (-1)^{r-1}}{\ds 
    (r-1)!}\int_{0}^{1}\log^{r-1}(x)\frac{\ds 1-x^n}{\ds 1-x}dx,
    \]
    and thus
    \[
    \begin{array}{lll}
	S_{r,q} & = & \frac{\ds (-1)^{r-1}}{\ds 
    (r-1)!}{\ds \int_{0}^{1}}\frac{\ds \log^{r-1}(x)}{\ds 
    1-x}{\ds \sum_{n=1}^{\infty}}\frac{\ds 1-x^{n}}{\ds n^q}dx\eqskip\\
     & = & \frac{\ds (-1)^{r-1}}{\ds 
    (r-1)!}{\ds \int_{0}^{1}}\frac{\ds \log^{r-1}(x)}{\ds 
    1-x}\left[\zeta(q)-\pol{q}(x)\right]dx\eqskip\\
    & = & \zeta(r)\zeta(q)-\frac{\ds (-1)^{r-1}}{\ds 
    (r-1)!}{\ds \int_{0}^{1}}\frac{\ds \log^{r-1}(x)}{\ds 1-x}\pol{q}(x)dx,
    \end{array}
    \]
    where the last step follows from (see identity 4.271
    in~\cite{grry}, for instance)
    \begin{equation}\label{logint}
    \int_{0}^{1}\frac{\ds \log^{r-1}(x)}{\ds 
    1-x} dx = (-1)^{r-1}\zeta(r)\Gamma(r).
    \end{equation}
\end{proof}

For the higher--order Euler sums $S_{1^{p},q}$ we have the following 
relation.
\begin{lem}\label{transf1}
    For $p\geq 1$ and $q\geq 2$ we have that
    \[
    S_{1^{p},q}= (-1)^{p}\int_{Q_{p}}\frac{\ds \pol{q-p}(x_{1}\ldots 
    x_{p})\prod_{j=1}^{p} \log(1-x_{j})}{\ds x_{1}\ldots x_{p}}dx,
    \]
    where $Q_{p}=(0,1)^{p}$.
\end{lem}
\begin{proof}
    Since $H_{n}=-n\dint_{0}^{1}y^{n-1}\log(1-y)dy$,
    we have that
    \[
    \begin{array}{lll}
    S_{1^{p},q} & = & \dsum_{n=1}^{\infty}\frac{\ds 
    (-1)^{p}n^{p}\left[\dint_{0}^{1}y^{n-1}\log(1-y)dy\right]^{p}}
    {\ds n^{q}}\eqskip\\
    & = & (-1)^{p}\dint_{Q_{p}}\dsum_{n=1}^{\infty}\frac{\ds 
    \left(x_{1}\ldots x_{p}\right)^{n-1}}{\ds 
    n^{q-p}}\log(1-x_{1})\ldots \log(1-x_{p})dx\eqskip\\
    & = & (-1)^{p}\dint_{Q_{p}}\frac{\ds 
    \pol{q-p}\left(x_{1}\ldots x_{p}\right)\log(1-x_{1})\ldots \log(1-x_{p})}
    {\ds x_{1}\ldots x_{p}}dx
    \end{array}
    \]
    as desired.
\end{proof}

Although this result will be used mainly to evaluate the integrals 
based on the knowledge of the value of the series, as was mentioned 
in the Introduction, in the Appendix we 
show how this can be yet another way of evaluating sums such as 
$S_{1^{2},2}$.

\section{\label{redsec}Reduction of integrals involving polylogarithms}

\subsection{Integrals related directly to $S_{r^{p},q}$}

We begin by presenting a simple example which is a direct 
consequence of Lemma~\ref{transf2} and the result in~\cite{bbg} for 
the case of $S_{r,q}$.

\begin{thm} \label{r0q} For even weight $w = q+r$, $(q\geq2, r\geq1)$ the integral
    \[
    K(r,0,q) = \int_{0}^{1}\frac{\ds \log^{r}(x)\pol{q}(x)}{\ds 1-x}dx
    \]
    is reducible to zeta values. More precisely,
    \[
    \begin{array}{lll}
    K(r,0,q) & = &
    r!\left\{ (-1)^{r+1}\zeta(w+1)\left[\frac{\ds 1}{\ds 
    2}-\frac{\ds (-1)^{r+1}}{\ds 2}\left(\begin{array}{c}w\\ 
    r+1\end{array}\right)-\frac{\ds (-1)^{r+1}}{\ds 2}\left(\begin{array}{c}w\\ 
    q\end{array}\right)\right]\right. \eqskip\\
     & & +\frac{\ds (-1)^{r}-1}{\ds 2}\zeta(r+1)\zeta(q)\eqskip\\
     & & +\dsum_{k=1}^{\lfloor (r+1)/2\rfloor}\left(\begin{array}{c}
     w-2k\\
     q-1\end{array}\right)\zeta(2k)\zeta(w-2k+1)\eqskip\\
     & & \left. +\dsum_{k=1}^{\lfloor q/2\rfloor}\left(\begin{array}{c} w-2k\\
     r\end{array}\right)\zeta(2k)\zeta(w-2k+1)\right\},
    \end{array}
    \]
    where $\zeta(1)$ should be interpreted as $0$ whenever it occurs.
\end{thm}
\begin{proof}
    This follows directly from Lemma~\ref{transf2} and the fact that 
    $S_{r,q}$ is reducible to zeta values for $q+r$ even~\cite{bbg} -- here
    we used the expression for
    $S_{r,q}$ given in~\cite{flsa}.
\end{proof}

When $w$ is odd it will no longer be possible in general to reduce
the integrals above to zeta values. This is in line with the discussion on
page 17 of~\cite{flsa}, for instance, which states that there are
{\it general} and {\it exceptional} classes of evaluations, in the 
sense that although there are families of sums which are not 
reducible to zeta values, some low weights are still reducible.
Thus, in the case of low odd
weights we may still obtain some explicit expressions for the above 
integrals such
as the following
\[
\begin{array}{l}
\dint_{0}^{1}\frac{\ds \log(x)\pol{2}(x)}{\ds 1-x}dx = -\frac{\ds 
3}{\ds 4}\zeta(4)\eqskip\\
\dint_{0}^{1}\frac{\ds \log^{2}(x)\pol{3}(x)}{\ds 1-x}dx = 
\zeta^2(3)-\zeta(6)\eqskip\\
\dint_{0}^{1}\frac{\ds \log(x)\pol{4}(x)}{\ds 1-x}dx = 
\zeta^2(3)-\frac{\ds 25}{\ds 12}\zeta(6)
\end{array}
\]
These may be derived from the expressions for the linear sums 
$S_{2,2}$, $S_{3,3}$ and $S_{2,4}$ given in~\cite{flsa}.

For other odd weights, and as is mentioned in~\cite{flsa}, it is possible
to obtain linear relations betwen different integrals. Here we give one example
obtained using the corresponding relation on page 23 
of~\cite{flsa} -- see also Table 9.
\[
\dint_{0}^{1}\frac{\ds \log^{2}(x)\pol{5}(x)}{\ds 1-x}dx -
5\dint_{0}^{1}\frac{\ds \log(x)\pol{6}(x)}{\ds 1-x}dx =
\frac{\ds 163}{\ds 12}\zeta(8)-8\zeta(3)\zeta(5).
\]

Another example is obtained by taking $p$ equal to $2$ in Lemma~\ref{transf1},
leading to a double integral which can then be reduced to a quadratic Euler sum 
and thus to linear sums and zeta values. This last step follows from 
the result in~\cite{bbg} regarding the reducibility of the sum 
$S_{1^{2},q}$ to linear sums and polynomials in zeta values for all 
weights. In the formulation of the result below, we used the 
expression for $S_{1^{2},q}$ given by Theorem 4.1 in~\cite{flsa} -- 
note that there is a misprint in the statement of Theorem 4.1 in~\cite{flsa}
where the residue term is missing.

\begin{thm} \label{mult}For any $q\geq0$ we have that
    \[
    \begin{array}{lll}
    \dint_{Q_{2}} \frac{\ds \pol{q}(xy)\log(1-x)\log(1-y)}{\ds 
    xy}dxdy & = & S_{2,q+2}+(q+2)S_{1,q+3}\eqskip\\
    & & \hspace*{.3cm}-\frac{\ds (q+2)(q+3)}{\ds 
    6}\zeta(q+2)\eqskip\\
    & & \hspace*{.6cm}+\zeta(2)\zeta(q+2)-\frac{\ds 1}{\ds 3}R(q),
    \end{array}
    \]
    where
    \[
    R(q) = \mbox{\rm Res}_{z=0}\left[\frac{\ds 
    \left(\psi(-z)+\gamma\right)^{3}}{\ds z^{q+2}}\right].
    \]
    Here $\psi$ is the logarithmic derivative of the Gamma function, 
    and $\gamma$ is Euler's constant.
\end{thm}

This gives, for instance,
    \[
    \begin{array}{lll}
    \dint_{Q_{2}} \frac{\ds \pol{3}(xy)\log(1-x)\log(1-y)}{\ds 
    xy}dxdy & = & -\frac{\ds 5}{\ds 
    2}\zeta(3)\zeta(4)-\zeta(2)\zeta(5)+6\zeta(7).
    \end{array}
    \]

We shall now consider the case of $J(-1,p,q)$, which may also be 
related directly to an Euler sum, and which will be used in the proof 
of Theorem~\ref{maint} for the remaining values of $m$.
\begin{thm} \label{jm1pq} For $p, q\geq 1$ the integral
    $
    J(-1,p,q) 
    $
    is reducible to zeta values. More precisely, for $p\geq q$, we have
    \[
    \begin{array}{l}
    J(-1,p,q) = (-1)^{q+1}(1+\frac{\ds p+q}{\ds 2})\zeta(p+q+1)+
    2(-1)^{q}\dsum_{j=1}^{\lfloor q/2\rfloor}\zeta(2j)\zeta(p+q-2j+1)\eqskip\\
    \hspace*{4cm}+
    \frac{\ds (-1)^{q}}{\ds 
    2}\dsum_{j=1}^{p-q}\zeta(j+q)\zeta(p-j+1).
    \end{array}
    \]
    A similar expression can be found for $p< q$.
\end{thm}
\begin{proof} Since $J(-1,p,q)=J(-1,q,p)$, we may, without loss of 
generality, assume that $p\geq q$.
    Writing $m=p-q\geq0$ and integrating by parts $q-1$ times gives that
    \[
    \begin{array}{lll}
    \dint_{0}^{1}\frac{\ds \pol{q+m}(x)\pol{q}(x)}{\ds x}dx & = &
    \dsum_{j=1}^{q-1}(-1)^{(j+1)}\zeta(q+m+j)\zeta(q-j+1)\eqskip\\
    & & \hspace*{1cm}+ (-1)^{q+1}
    \dint_{0}^{1}\frac{\ds \pol{2q+m-1}(x)\pol{1}(x)}{\ds 
    x}dx\eqskip\\
    & = & \dsum_{j=1}^{q-1}(-1)^{(j+1)}\zeta(q+m+j)\zeta(q-j+1)
    \eqskip\\
    & & \hspace*{1cm}+
    (-1)^{q+1}S_{1,2q+m},
    \end{array}
    \]
    where the last step follows from Lemma~\ref{transf1}. Using the 
    expression for $S_{1,2q+m}$ (already known to Euler), we get
    \[
    \begin{array}{lll}
    \dint_{0}^{1}\frac{\ds \pol{q+m}(x)\pol{q}(x)}{\ds x}dx & = &
    \dsum_{j=1}^{q-1}(-1)^{(j+1)}\zeta(q+m+j)\zeta(q-j+1)\eqskip\\
    & & \hspace*{1cm}+(-1)^{q+1}\left[ (1+\frac{\ds 2q+m}{\ds 
    2})\zeta(2q+m+1)\right.    \eqskip\\
    & & \hspace*{1cm}-\left.\frac{\ds 1}{\ds 
    2}\dsum_{j=1}^{2q+m-2}\zeta(j+1)\zeta(2q+m-j)\right]\eqskip\\
    & = & (-1)^{q+1}(q+\frac{\ds m}{\ds 
    2}+1)\zeta(2q+m+1)\eqskip\\
    & & \hspace*{1cm}+2(-1)^{q}\dsum_{j=1}^{\lfloor 
    q/2\rfloor}\zeta(2j)\zeta(2(q-j)+m+1) \eqskip\\
    & & \hspace*{1cm}+
    \frac{\ds (-1)^{q}}{\ds 2}\dsum_{j=1}^{m}\zeta(j+q)\zeta(q+m-j+1).
    \end{array}
    \]
\end{proof}

Letting $p=q$ gives the following
\begin{cor}\label{cor2}
    \[
    \int_{0}^{1} \frac{\pol{p}^{2}(x)}{x}dx = 
    (-1)^{p+1}(p+1)\zeta(2p+1)+2(-1)^{p}\dsum_{j=1}^{\lfloor 
    p/2\rfloor} \zeta(2j)\zeta[2(p-j)+1].
    \]
\end{cor}

\subsection{Recurrence relations}

We shall now develop the recurrence relations that will allow us to 
reduce the integrals $J(m,p,q)$ and $K(r,p,q)$.
\begin{lem}\label{reclem}
    For $p,q\geq2$ and $m\geq 0$, the following recurrence relations hold:
    \begin{enumerate}
	\item[(i)] $
	J(m,p,q) = \frac{\ds \zeta(p)\zeta(q)}{\ds m+1}-\frac{\ds 1}
	{\ds m+1}\left[J(m,p-1,q)+J(m,p,q-1)\right]\eqskip
	$
	\item[(ii)]
	$
	\left.
	\begin{array}{lll}
	J(m,1,q) & = & \frac{\ds \zeta(q)}{\ds m+1}-\frac{\ds 1}{\ds m+1}
	\left[ m J_{0}(m,q) + J_{0}(m,q-1)\right]\eqskip\\
	& & \hspace*{.5cm}+
	\frac{\ds 1}{\ds m+1}\left[ m J(m-1,1,q)+J(m-1,1,q-1)\right.
	\eqskip\\
	& & \hspace*{1cm}\left. -J(m,1,q-1)
	\right],
	\end{array}
	\right\}\eqskip
	$
	
	where $J_{0}(m,q) := \dint_{0}^{1}x^{m}\pol{q}(x)dx.\eqskip$
	\item[(iii)]
	$
	J_{0}(m,q) = \frac{\ds \zeta(q)}{\ds m+1}-\frac{\ds 1}{\ds m+1}
	J_{0}(m,q-1).
	$
    \end{enumerate}
\end{lem}
\begin{proof} The proof of the three relations follows by integration 
by parts, in the first and third cases using the integral of $x^{m}$, 
and in the second case the integral of $\log(1-x)$.
\end{proof}

\begin{rem} {\rm Note that the first relation in the above Lemma also 
holds in the case of $m=-2$. Also, in the caso where $m=0$, (ii)
simplifies to
    \[
    J(0,1,q) = \zeta(q) + \dsum_{k=2}^{q-1}(-1)^{q+k}\zeta(k)- 
    (-1)^{q}-J(0,1,q-1) + J(-1,1,q-1),
    \]
where we have used the fact that
    \[
    \dint_{0}^{1} \pol{q}(x)dx = 
    (-1)^{q-1}+\dsum_{k=2}^{q}(-1)^{k+q}\zeta(k).
    \]
}
\end{rem}

Successive application of these relations in the case of 
non--negative $m$ leads to expressions 
containing only the integrals 
$J(-1,1,q), J(m,1,1)$ and $J_{0}(m,1)$ -- note that $J(0,p,q)=J(0,q,p)$.
For $m=-2$, we cannot apply the second recurrence in the Lemma, and it
becomes necessary to be able to evaluate $J(-2,1,q)$ separately.
Thus, in order to be possible to actually use the relations in Lemma~\ref{reclem}
to evaluate $J(m,p,q)$, we have to determine the values of 
the above four integrals. 

The integral $J(-1,1,q)$ is given by Theorem~\ref{jm1pq}.

In the case of $J_{0}(m,1)$ we have from the proof of 
Lemma~\ref{transf1} that
\[
-\dint_{0}^{1} x^{m}\log(1-x)dx = \frac{\ds 1}{\ds m+1}H_{m+1},
\]
while for $J(m,1,1)$ identity (4.2.5) in~\cite{dedu} gives that
\[
\dint_{0}^{1}x^{m} \log^{2}(1-x)dx = \frac{\ds 2}{\ds m+1}\left[
H_{m+1}^{(2)}+\dsum_{k=1}^{m+1}\frac{\ds H_{k}}{\ds k+1}\right].
\]

It remains to determine $J(-2,1,q)$. This can again be done
by means of a recurrence relation, which relates these integrals to $J(-1,1,q-1)$.
\begin{lem}\label{reclem2}
    For $q\geq 2$ the following recurrence relation holds:
    \[
    J(-2,1,q) = \zeta(q+1)+J(-2,1,q-1)-J(-1,1,q-1).\eqskip
    \]
    Furthermore, $J(-2,1,1)=2\zeta(2)$.
\end{lem}
\begin{proof}
    One again proceeds by integration by parts, but now starting from 
    $J(-1,1,q)$ and using the integral of 
    $\log(1-x)$. The terms in $J(-1,1,q)$ cancel out, and the result 
    follows. The value of $J(-2,1,1)$ can be found in~\cite{dedu}, for 
    instance.
\end{proof}

For the case of $K$ we need two recurrence relations. The first is 
similar to those used in the previous case, and is again based on 
integration by parts.
\begin{lem}\label{krecur}
    For positive integers $p,q$ and $r$ we 
    have
    \[
    K(r,p,q) = -\frac{\ds 1}{\ds r+1}\left[ 
    K(r+1,p-1,q)+K(r+1,p,q-1)\right].
    \]
\end{lem}
\begin{proof} Integrate by parts, using the integral of 
$\log^{r}(x)/x$.
\end{proof}

The second recursion is based on the following symmetry relation for 
linear sums
\begin{equation}\label{symmrelser}
S_{r,q}+S_{q,r}=\zeta(r)\zeta(q) + \zeta(r+q).
\end{equation}

\begin{lem}\label{symmrec}
    For $r\geq 1,q\geq2$ we have that
    \[
    K(r,0,q) = (-1)^{r+q}\frac{\ds r!}{\ds 
    (q-1)!}K(q-1,0,r+1)+(-1)^{r}r!\left[ 
    \zeta(r+1)\zeta(q)-\zeta(r+q+1)\right].
    \]
    In particular.
    \[
    K(r,0,r+1) = \dint_{0}^{1}\frac{\ds \log^{r}(x)\pol{r+1}(x)}{\ds 
    1-x}dx =  (-1)^{r+1}\frac{\ds r!}{\ds 2}\left[ 
    \zeta(2r+2)-\zeta^{2}(r+1)\right].
    \]
\end{lem}
\begin{proof}
    The first expression is just a {\it translation} of the symmetry 
    relations~(\ref{symmrelser}) into the case of integrals using 
    Lemma~\ref{transf2}. Applying this expression to the case 
    of $q=r+1$ gives the second statement.
\end{proof}

\subsection{\label{pt12}Proof of Theorems~\ref{maint} and~\ref{maint2}}

By combining the above results, it is now possible to prove Theorems~\ref{maint}
and~\ref{maint2}.

{\it Proof of Theorem~\ref{maint}.} As was mentioned above, the 
application of the recurrence relations in Lemma~\ref{reclem} allows us 
to reduce $J(m,p,q)$ to rational combinations of zeta values and
integrals of the type $J(m,1,1), J_{0}(m,1)$ and $J(0,1,q)$. The 
first two of these reduce to rational numbers (see identities (4.2.4) 
and (4.2.5) in~\cite{dedu}, for instance), while the third can, by 
Lemma~\ref{reclem2}, be written in 
terms of an integer, zeta values, and $J(-1,1,q)$, which is, in turn, 
reducible to zeta values by Theorem~\ref{jm1pq}.
\qed

{\it Proof of Theorem~\ref{maint2}.} By Lemma~\ref{krecur} we have that the
weight $w=p+q+r$ does not change as we successivly apply the 
recurrence relation. This, together with the fact that $K(r,p,q) = 
K(r,q,p)$, allows us to reduce the integral $K(r,p,q)$ to integrals of 
the form $K(r',0,q')$, with $r'+q' = p+q+r$. By 
Theorem~\ref{r0q}, we have that these integrals are reducible if 
$r'+q'$ is even and $q'\geq2$. On the other hand, the case $K(r,0,1)$ 
is reducible for all $r$ (see~\cite{kolb}, for instance, and the 
references therein), and so the result
follows.

In the case of odd weights, by applying directly the results 
from~\cite{bbg,zag} we get the $\lfloor (w-1)/6\rfloor$ bound. 

We now use the following version of 
Lemma~\ref{krecur}
\[
K(r,p,q) = -r K(r-1,p,q+1)-K(r,p-1,q+1).
\]
By applying this recursion repeatedly we see that it is possible to reduce
the original integral 
to integrals of the form $K(0,p',q')=J(-1,p',q')$ (reducible by 
Theorem~\ref{maint}), and integrals $K(i,0,p+q+r-i)$ with 
$i=1,\ldots,r$.

We shall now prove that when $i$ is even, $K(i,0,q')$ can be 
reduced to integrals where the first entry is $i-1$. To see this, note that 
for even $i$ we have $K(i,0,q')$ with $q'$ odd. If we apply the 
above recurrence to $K(i,q',0)$ repeatedly, we get a linear combination of 
integrals with the first entry $i-1$, plus one of the form 
$K(i,p'+1,p')=K(i,p',p'+1)$. Aplying the recursion once more to this 
(or by letting $p=q=p'+1$ in Lemma~\ref{krecur}) we get
\[
K(i,p',p'+1) = -\frac{\ds i}{\ds 2}K(i-1,p'+1,p'+1)
\]
which enables us to reduce this last term to one with odd first entry
as desired.

Finally, the relations in Lemma~\ref{symmrec} allow us to reduce any 
$K(r,0,q)$ with $q>r$ to one with $q<r-2$.
\qed

\subsection{\label{pt3}Estimates for the integrals $K(r,0,q)$}

We begin by estimating the Euler sums $S_{rq}$, and then obtain 
the result for $\kappa_{rq}$ from these. We have
\[
\begin{array}{lll}
    S_{rq} & = & \dsum_{n=1}^{\infty}\frac{\ds 1}{\ds n^{q}}
    \left(1+\frac{\ds 1}{\ds 2^{r}}+\ldots+\frac{\ds 1}{\ds 
    n^{r}}\right)\eqskip\\
    & \leq & \zeta(q)+\frac{\ds 1}{\ds 2^{r}}\dsum_{n=1}^{\ds 
    \infty}\frac{\ds n-1}{\ds n^{q}}\eqskip\\
    & \leq &\zeta(q)+\frac{\ds 1}{\ds 2^{r}}\left[\zeta(q-1)-\zeta(q)\right],
\end{array}
\]
and so
\[
0\leq S_{rq}-\zeta(q)\leq \frac{\ds 1}{\ds 2^{r}}\left[\zeta(q-1)-
\zeta(q)\right].
\]
Theorem~\ref{maint3} now follows by combining this with Lemma~\ref{transf2}.

\section{\label{examp}Some examples}

We begin by considering the integrals $J(m,p,q)$, for $-2\leq m\leq 1$ 
and $2\leq p,q\leq 4$.
    \[
    \begin{array}{rcl}
    \dint_{0}^{1} \frac{\ds \pol{2}^{2}(x)}{\ds x^{2}}dx & = & 
    4\zeta(2)-2\zeta(3)-\frac{\ds 5}{\ds 2}\zeta(4)\eqskip\\
    \dint_{0}^{1} \frac{\ds \pol{2}(x)\pol{3}(x)}{\ds x^{2}}dx & = & 
    6\zeta(2)-3\zeta(3)-\frac{\ds 11}{\ds 4}\zeta(4)
    -\zeta(2)\zeta(3)\eqskip\\
    \dint_{0}^{1} \frac{\ds \pol{2}(x)\pol{4}(x)}{\ds x^{2}}dx & = & 
    8\zeta(2)-4\zeta(3)-3\zeta(4)
    -2\zeta(5)-\frac{\ds 7}{\ds 4}\zeta(6)\eqskip\\
    \dint_{0}^{1} \frac{\ds \pol{3}^{2}(x)}{\ds x^{2}}dx & = & 
    12\zeta(2)-6\zeta(3)-\frac{\ds 11}{\ds 2} 
    \zeta(4)-2\zeta(2)\zeta(3)-\zeta^{2}(3)\eqskip\\
	\dint_{0}^{1} \frac{\ds \pol{3}(x)\pol{4}(x)}{\ds x^{2}}dx & = &
	20\zeta(2)-10\zeta(3)-\frac{\ds 17}{\ds 2}\zeta(4)
	-2\zeta(5)-2\zeta(2)\zeta(3)\eqskip\\
	& & \hspace*{1cm} -\frac{\ds 7}{\ds 4}\zeta(6)
	-\zeta^{2}(3)-\zeta(3)\zeta(4)
    \eqskip\\
	\dint_{0}^{1} \frac{\ds \pol{4}^{2}(x)}{\ds x^{2}}dx & = &
	40\zeta(2)-20\zeta(3)-17\zeta(4)
	-4\zeta(5)-4\zeta(2)\zeta(3)\eqskip\\
	& & \hspace*{1cm} -\frac{\ds 7}{\ds 2}\zeta(6)
	-2\zeta^{2}(3)-2\zeta(3)\zeta(4)-\frac{\ds 7}{\ds 6}\zeta(8)
    \end{array}
    \]
    \begin{center}
	Table 1. $J(-2,p,q)$, $2\leq p,q\leq 4.$
    \end{center}
    \[
    \begin{array}{rcl}
    \dint_{0}^{1} \frac{\ds \pol{2}^{2}(x)}{\ds x}dx & = & 
    2\zeta(2)\zeta(3)-3\zeta(5)\eqskip\\
    \dint_{0}^{1} \frac{\ds \pol{2}(x)\pol{3}(x)}{\ds x}dx & = & 
    \frac{\ds 1}{\ds 2} \zeta^{2}(3)\eqskip\\
    \dint_{0}^{1} \frac{\ds \pol{2}(x)\pol{4}(x)}{\ds x}dx & = & 
    2\zeta(2)\zeta(5)+\zeta(3)\zeta(4)-4\zeta(7)\eqskip\\
    \dint_{0}^{1} \frac{\ds \pol{3}^{2}(x)}{\ds x}dx & = & 
    -2\zeta(2)\zeta(5)+4\zeta(7)\eqskip\\
    \dint_{0}^{1} \frac{\ds \pol{3}(x)\pol{4}(x)}{\ds x}dx & = & 
    \frac{\ds 7}{\ds 12} \zeta(8)\eqskip\\
    \dint_{0}^{1} \frac{\ds \pol{4}^{2}(x)}{\ds x}dx & = & 
    2\zeta(4)\zeta(5)+2\zeta(2)\zeta(7)-5\zeta(9)
    \end{array}
    \]
    \begin{center}
	Table 2. $J(-1,p,q)$, $2\leq p,q\leq 4.$
    \end{center}
    \[
    \begin{array}{rcl}
    \dint_{0}^{1} \pol{2}^{2}(x)dx & = &
    6-2\zeta(2)-4\zeta(3)+\frac{\ds 5}{\ds 2}\zeta(4)\eqskip\\
    \dint_{0}^{1} \pol{2}(x)\pol{3}(x)dx & = & -10 + 
    4\zeta(2)+5\zeta(3)-\frac{\ds 15}{\ds 4} \zeta(4)+\zeta(2)\zeta(3)
    \eqskip\\
    \dint_{0}^{1} \pol{2}(x)\pol{4}(x)dx & = & 15 - 
    7\zeta(2)-5\zeta(3)+4\zeta(4)-3\zeta(5)+\frac{\ds 7}{\ds 4} 
    \zeta(6)
    \eqskip\\
    \dint_{0}^{1} \pol{3}^{2}(x)dx & = & 20-8\zeta(2)-10\zeta(3)-
    \frac{\ds 15}{\ds 2} \zeta(4)-2\zeta(2)\zeta(3)+\zeta^{2}(3)
    \eqskip\\
    \dint_{0}^{1} \pol{3}(x)\pol{4}(x)dx & = & 
    -35+15\zeta(2)+15\zeta(3)-\frac{\ds 23}{\ds 2}\zeta(4)+3\zeta(5)+
    2\zeta(2)\zeta(3)\eqskip\\
    & & \hspace*{1cm} -\zeta^{2}(3)-\frac{\ds 7}{\ds 4}\zeta(6)+
    \zeta(3)\zeta(4)\eqskip\\
    \dint_{0}^{1} \pol{4}^{2}(x)dx & = & 70-30\zeta(2)-30\zeta(3)+23
    \zeta(4)-6\zeta(5)-4\zeta(2)\zeta(3)\eqskip\\
    & & \hspace*{1cm} +\frac{\ds 7}{\ds 2}\zeta(6)+ 2\zeta^{2}(3)
    -2\zeta(3)\zeta(4)+\frac{\ds 7}{\ds 
    6}\zeta(8)
    \end{array}
    \]
    \begin{center}
	Table 3. $J(0,p,q)$, $2\leq p,q\leq 4.$
    \end{center}
    \[
    \begin{array}{rcl}
    \dint_{0}^{1} x\pol{2}^{2}(x) dx & = & \frac{\ds 25}{\ds 16}-\frac{\ds 
    3}{\ds 4}\zeta(2)-\zeta(3)+\frac{\ds 5}{\ds 4}\zeta(4)\eqskip\\
    \dint_{0}^{1} x\pol{2}(x)\pol{3}(x) dx & = & -\frac{\ds 47}{\ds 32}+\frac{\ds 
    7}{\ds 8}\zeta(2)+\frac{\ds 3}{\ds 8}\zeta(3)-\frac{\ds 15}{\ds 16}\zeta(4)
    +\frac{\ds 1}{\ds 2}\zeta(2)\zeta(3)\eqskip\\
    \dint_{0}^{1} x\pol{2}(x)\pol{4}(x) dx & = & \frac{\ds 173}{\ds 128}-\frac{\ds 
    31}{\ds 32}\zeta(2)+\frac{\ds 3}{\ds 16}\zeta(3)+\frac{\ds 1}{\ds 4}\zeta(4)
    -\frac{\ds 3}{\ds 4}\zeta(5)+\frac{\ds 7}{\ds 8}\zeta(6)\eqskip\\
    \dint_{0}^{1} x\pol{3}^{2}(x) dx & = & \frac{\ds 47}{\ds 32}-\frac{\ds 
    7}{\ds 8}\zeta(2)-\frac{\ds 3}{\ds 8}\zeta(3)+\frac{\ds 15}{\ds 16}\zeta(4)
    -\frac{\ds 1}{\ds 2}\zeta(2)\zeta(3)+\frac{\ds 1}{\ds 
    2}\zeta^{2}(3)\eqskip\\
    \dint_{0}^{1} x\pol{3}(x)\pol{4}(x)dx & = & -\frac{\ds 361}{\ds 256}
    +\frac{\ds 59}{\ds 64}\zeta(2)+\frac{\ds 3}{\ds 32}\zeta(3)-
    \frac{\ds 19}{\ds 32}
    \zeta(4)+\frac{\ds 1}{\ds 4}\zeta(2)\zeta(3)
    \eqskip\\
    & & \hspace*{1cm} +\frac{\ds 3}{\ds 8}\zeta(5)-\frac{\ds 7}{\ds 16}\zeta(6)-\frac{\ds 1}{\ds 4}
    \zeta^{2}(3)+\frac{\ds 1}{\ds 2}\zeta(3)\zeta(4)
    \eqskip\\
    \dint_{0}^{1} x\pol{4}^{2}(x)dx & = & \frac{\ds 361}{\ds 256}
    -\frac{\ds 59}{\ds 64}\zeta(2)-\frac{\ds 3}{\ds 32}\zeta(3)+
    \frac{\ds 19}{\ds 32}
    \zeta(4)-\frac{\ds 1}{\ds 4}\zeta(2)\zeta(3)
    \eqskip\\
    & & \hspace*{1cm} -\frac{\ds 3}{\ds 8}\zeta(5)+\frac{\ds 7}{\ds 16}\zeta(6)
    +\frac{\ds 1}{\ds 4}\zeta^{2}(3)-\frac{\ds 1}{\ds 2}\zeta(3)\zeta(4)+
    \frac{\ds 7}{\ds 12}\zeta(8)
    \end{array}
    \]
    \begin{center}
	Table 4. $J(1,p,q)$, $2\leq p,q\leq 4.$
    \end{center}

We now give some values of the integrals $K(r,p,q)$ for $3\leq 
p+q+r\leq 7$\vspace*{.5cm}.

     \[
    \begin{array}{rcl}
	\dint_{0}^{1} \frac{\ds \log(x)\log^2(1-x)}{\ds x}dx & = & 
	-\frac{\ds 1}{\ds 2}\zeta(4)\eqskip\\
	\dint_{0}^{1} \frac{\ds \log(x)\pol{2}(x)}{\ds 1-x}dx & = & 
	-\frac{\ds 3}{\ds 4}\zeta(4)\eqskip\\
	\dint_{0}^{1} \frac{\ds \log^{2}(x)\log(1-x)}{\ds 1-x}dx & = & 
	-\frac{\ds 1}{\ds 2}\zeta(4)\eqskip\\
    \end{array}
    \]
    \begin{center}
	Table 5. $K(r,p,q)$, $p+q+r=3$.
    \end{center}
    \[
    \begin{array}{rcl}
	\dint_{0}^{1} \frac{\ds \log(x)\pol{3}(x)}{\ds 1-x}dx & = & 
	2\zeta(2)\zeta(3)-\frac{\ds 9}{\ds 2}\zeta(5)\eqskip\\
	\dint_{0}^{1} \frac{\ds \log(x)\log(1-x)\pol{2}(x)}{\ds x}dx & = & 
	\zeta(2)\zeta(3)-\frac{\ds 9}{\ds 6}\zeta(5)\eqskip\\
	\dint_{0}^{1} \frac{\ds \log^{2}(x)\pol{2}(x)}{\ds 1-x}dx & = & 
	6\zeta(2)\zeta(3)-11\zeta(5)\eqskip\\	
	\dint_{0}^{1} \frac{\ds \log^{2}(x)\log^2(1-x)}{\ds x}dx & = & 
	-4\zeta(2)\zeta(3)+8\zeta(5)\eqskip\\
	\dint_{0}^{1} \frac{\ds \log^{3}(x)\log(1-x)}{\ds 1-x}dx & = & 
	-6\zeta(2)\zeta(3)+12\zeta(5)\eqskip\\
    \end{array}
    \]
    \begin{center}
	Table 6. $K(r,p,q)$, $p+q+r=4$.
    \end{center}
    
    \[
    \begin{array}{rcl}
	\dint_{0}^{1} \frac{\ds \log(x)\pol{4}(x)}{\ds 1-x}dx & = & 
	-\frac{\ds 25}{\ds 12}\zeta(6)+\zeta^{2}(3)\eqskip\\
	\dint_{0}^{1} \frac{\ds \log(x)\log(1-x)\pol{3}(x)}{\ds x}dx & = & 
	-\frac{\ds 1}{\ds 3}\zeta(6)+\frac{\ds 1}{\ds 2}\zeta^{2}(3)\eqskip\\
	\dint_{0}^{1} \frac{\ds \log(x)\pol{2}^{2}(x)}{\ds x}dx & = & 
	-\frac{\ds 1}{\ds 3}\zeta(6)\eqskip\\
	\dint_{0}^{1} \frac{\ds \log^{2}(x)\pol{3}(x)}{\ds 1-x}dx & = & 
	-\zeta(6)+\zeta^{2}(3)\eqskip\\	
	\dint_{0}^{1} \frac{\ds \log^{2}(x)\log(1-x)\pol{2}(x)}{\ds x}dx & = & 
	-\frac{\ds 1}{\ds 3}\zeta(6)\eqskip\\
	\dint_{0}^{1} \frac{\ds \log^{3}(x)\pol{2}(x)}{\ds 1-x}dx & = & 
	8\zeta(6)-6\zeta^{2}(3)\eqskip\\
	\dint_{0}^{1} \frac{\ds \log^{3}(x)\log^2(1-x)}{\ds x}dx & = & 
	-9\zeta(6)+6\zeta^{2}(3)\eqskip\\
	\dint_{0}^{1} \frac{\ds \log^{4}(x)\log(1-x)}{\ds 1-x}dx & = & 
	-18\zeta(6)+12\zeta^{2}(3)\eqskip\\
    \end{array}
    \]
    \begin{center}
	Table 7. $K(r,p,q)$, $p+q+r=5$.
    \end{center}
    
    \[
    \begin{array}{rcl}
	\dint_{0}^{1} \frac{\ds \log(x)\pol{5}(x)}{\ds 1-x}dx & = & 
	2\zeta(3)\zeta(4)+4\zeta(2)\zeta(5)-10\zeta(7)\eqskip\\
	\dint_{0}^{1} \frac{\ds \log(x)\log(1-x)\pol{4}(x)}{\ds x}dx & = & 
	\zeta(3)\zeta(4)+3\zeta(2)\zeta(5)-6\zeta(7)\eqskip\\
	\dint_{0}^{1} \frac{\ds \log(x)\pol{2}(x)\pol{3}(x)}{\ds x}dx & = & 
	\zeta(2)\zeta(5)-2\zeta(7)\eqskip\\
	\dint_{0}^{1} \frac{\ds \log^{2}(x)\pol{4}(x)}{\ds 1-x}dx & = & 
	2\zeta(3)\zeta(4)+20\zeta(2)\zeta(5)-36\zeta(7)\eqskip\\
	\dint_{0}^{1} \frac{\ds \log^{2}(x)\log(1-x)\pol{3}(x)}{\ds x}dx & = & 
	14\zeta(2)\zeta(5)-24\zeta(7)\eqskip\\
	\dint_{0}^{1} \frac{\ds \log^{2}(x)\pol{2}^{2}(x)}{\ds x}dx & = & 
	12\zeta(2)\zeta(5)-20\zeta(7)\eqskip\\
	\dint_{0}^{1} \frac{\ds \log^{3}(x)\pol{3}(x)}{\ds 1-x}dx & = & 
	60\zeta(2)\zeta(5)-102\zeta(7)\eqskip\\	
	\dint_{0}^{1} \frac{\ds \log^{3}(x)\log(1-x)\pol{2}(x)}{\ds x}dx & = & 
	-18\zeta(2)\zeta(5)+30\zeta(7)\eqskip\\
	\dint_{0}^{1} \frac{\ds \log^{4}(x)\pol{2}(x)}{\ds 1-x}dx & = & 
	120\zeta(2)\zeta(5)+48\zeta(3)\zeta(4)-264\zeta(7)\eqskip\\
	\dint_{0}^{1} \frac{\ds \log^{4}(x)\log^{2}(1-x)}{\ds x}dx & = & 
	-48\zeta(2)\zeta(5)-48\zeta(3)\zeta(4)+144\zeta(7)\eqskip\\	
	\dint_{0}^{1} \frac{\ds \log^{5}(x)\log(1-x)}{\ds 1-x}dx & = & 
	-120\zeta(2)\zeta(5)-120\zeta(3)\zeta(4)+360\zeta(7)\eqskip\\	
    \end{array}
    \]
    \begin{center}
	Table 8. $K(r,p,q)$, $p+q+r=6$.
    \end{center}	

    The next table presents the values of the integrals $K(r,0,q)$, 
    with $r+q=7$. In this case, all values can be reduced to zeta 
    values at integers plus a new constant, which we took to be 
    $\kappa_{16}= K(1,0,6)\approx -0.651565$.
    Theorem~\ref{maint3} yields
    \[
    \left|\kappa_{16}-\zeta(6)\left[1-\zeta(2)\right]\right|\leq
    \frac{\ds \zeta(5)-\zeta(6)}{\ds 4}\approx 0.0049.
    \]
    \[
    \begin{array}{rcl}
	\dint_{0}^{1} \frac{\ds \log(x)\pol{6}(x)}{\ds 1-x} & = & \kappa_{16}
	\eqskip\\
	\dint_{0}^{1} \frac{\ds \log^{2}(x)\pol{5}(x)}{\ds 1-x} & = &
	\frac{\ds 163}{\ds 12}\zeta(8)+5\kappa_{16}-8\zeta(3)\zeta(5)
	\eqskip\\
	\dint_{0}^{1} \frac{\ds \log^{3}(x)\pol{4}(x)}{\ds 1-x} & = &
	-\frac{\ds \zeta(8)}{\ds 2}
	\eqskip\\
	\dint_{0}^{1} \frac{\ds \log^{4}(x)\pol{3}(x)}{\ds 1-x} & = &
	-187\zeta(8)-60\kappa_{16}+120\zeta(3)\zeta(5)\eqskip\\
	\dint_{0}^{1} \frac{\ds \log^{5}(x)\pol{2}(x)}{\ds 1-x} & = &
	-80\zeta(8)-120\kappa_{16}\eqskip\\
	\dint_{0}^{1} \frac{\ds \log^{6}(x)\log(1-x)}{\ds 1-x} & = &
	720\zeta(3)\zeta(5)-900\zeta(8)\eqskip\\
    \end{array}
    \]
    \begin{center}
	Table 9. $K(r,0,q)$, $q+r=7$.
    \end{center}	
\appendix

\section{The Euler sum $S_{1^{2},2}$ and a related double integral}

The Euler sum
\[
\sum_{n=1}^{\infty}\frac{\ds H_{n}^{2}}{\ds (n+1)^2}
\]
was shown in~\cite{bobo} to be equal to the integral
\[
\frac{\ds 1}{\ds \pi}\dint_{0}^{\pi} t^2\log^{2}\left[2\cos(\frac{\ds 
t}{\ds 2})\right]dt,
\]
by using an appropriate Fourier series and then Parseval's identity
-- for other ways of evaluating this series, see~\cite{baap,chu,doel,flsa}. 
The above integral was then evaluated 
via contour integration -- see also Section~7.9 in~\cite{lewi},
specially~7.9.10 --, and used to establish the value of $S_{1^{2},2}$. Here we use 
Lemma~\ref{transf1} to transform the 
series $S_{1^{2},2}$ into a double integral which we then evaluate directly.

\begin{prop}
    \[
\dsum_{n=1}^{\infty} \left(\frac{\ds H_{n}}{\ds n}\right)^{2} =
\frac{\ds 17\zeta(4)}{\ds 4}.
\]
\end{prop}
\begin{proof}
By Lemma~\ref{transf1} we have that
\[
\dsum_{n=1}^{\infty} \left(\frac{\ds H_{n}}{\ds n}\right)^{2} = 
\dint_{0}^{1}\dint_{0}^{1} \frac{\ds \log(1-x)\log(1-y)}{\ds 1-xy}dxdy.
\]
On the other hand,
\[
\begin{array}{lll}
    \dint_{0}^{1}\dint_{0}^{1} \frac{\ds \log(1-x)\log(1-y)}{\ds 
    1-xy}dxdy & = &\dint_{0}^{1}\log(1-y)\dint_{0}^{1}\frac{\ds \log(1-x)}{\ds 
    1-xy}dxdy
\end{array}
\]
\[
\begin{array}{l}
    = -\dint_{0}^{1}\log(1-y)\left.\frac{\ds 
    \log(1-x)\log\left(\frac{\ds xy-1}{\ds 
    y-1}\right)+\pol2\left[\frac{\ds y(1-x)}{\ds y-1}\right]}{\ds 
    y}\right]^{1}_{0}dy\eqskip\\
    =\dint_{0}^{1}\frac{\ds \log(1-y)}{\ds y}\pol2\left(\frac{\ds 
    y}{\ds y-1}\right)dy.
\end{array}
\]
Since, from~\cite{lewi}, we have that
\[
\pol2\left(\frac{\ds y}{\ds y-1}\right) = - \frac{\ds 1}{\ds 
2}\log^{2}(1-y)-\pol2(y),
\]
the above integral is equal to
\[
\begin{array}{lll}
-\dint_{0}^{1}\frac{\ds \log(1-y)}{\ds y}\left[
\frac{\ds 1}{\ds 2}\log^{2}(1-y)+\pol2(y)\right]& = &
-\frac{\ds 1}{\ds 2}\dint_{0}^{1}\frac{\ds \log^{3}(1-y)}{\ds y}\eqskip\\
& & \ \ \ \ \  -\dint_{0}^{1}\frac{\ds \log(1-y)}{\ds y}\pol2(y)dy
\end{array}
\]
From~(\ref{logint}) it follows that the first of these integrals is
equal to $-6\zeta(4)$. In the case of the second integral, the 
integrand is the derivative of $-\pol2^{2}(y)/2$, and so we get 
$-10\zeta(4)/4$. Substituting this back we obtain the desired result.
\end{proof}

\end{document}